\documentclass[english,11pt,a4paper,leqno]{article}
\usepackage{babel}
\usepackage{bm}
\usepackage{makeidx}
\usepackage{layout}
\usepackage[latin1]{inputenc}
\usepackage{latexsym,amssymb,amsmath}
\usepackage{amsfonts,amsbsy}
\usepackage{eufrak}
\usepackage[mathscr]{eucal}
\usepackage{verbatim}
\usepackage{amscd}
\usepackage{xy}
\xyoption{all}
\usepackage{indentfirst}
\numberwithin{equation}{section}
\begin{document}
\date{}
\title{\bf {\large{RELATION OF THE CYCLOTOMIC EQUATION WITH THE HARMONIC AND  DERIVED  SERIES}}}
\author{LUIS J. BOYA \footnote{\texttt{Departamento de F\'{\i}sica Te\'{o}rica, Universidad de Zaragoza E-50009 Zaragoza, SPAIN(luis.boya@gmail.com)}}\quad AND \quad CRISTIAN RIVERA\footnote{\texttt{Departamento de F\'{\i}sica Te\'{o}rica, Universidad de Zaragoza E-50009 Zaragoza, SPAIN(cristian\_elfisico@hotmail.com)}} }
\maketitle


\textbf{Mathematical science classification.}   11 Y 40, 	  40 A 05,          41 A 55.\\


\begin{abstract}
We associate some (old) convergent series related to definite integrals with the cyclotomic equation $x^m-1= 0$, for several natural numbers $m$; for example, for $m = 3$, $x^3-1 = (x-1)(1+x+x^2)$, leads to $\int_0^1dx\frac{1}{(1+x+x^2)} = \frac{\pi}{(3\sqrt{3})} = (1-\frac{1}{2}) + (\frac{1}{4}-\frac{1}{5}) +       (\frac{1}{7}-\frac{1}{8}) + \ldots$  . In some cases, we express the results in terms of the Dirichlet characters. Generalizations for arbitrary $m$ are well defined, but do imply integrals and/or series summations rather involved.
\\
\end{abstract}

\section{Introduction} Nicola Oresmes proved ca. 1350, the \emph{divergence} of the \emph{harmonic series} (e.g., see [1], p. 183). Indeed, today we know precisely how it does \emph{diverge} ([2], p. 14; see also [3]):

\begin{equation}
\textrm{Harm(N)}:=\sum_{n=1}^N \frac{1}{n}
 \end{equation}

 For $N>>1$

 \begin{equation}
\textrm{Harm(N)}\longrightarrow log(N)+\gamma+\frac{1}{2N}-\frac{1}{12N(N+1)}-\ldots
 \end{equation}

\noindent where $\gamma$ is the Euler-Mascheroni constant (written and named $C$ at times; for a recent reference, see [4])

\begin{equation}
\gamma=\lim_{N\longrightarrow\infty}\left(\sum_{n=1}^{N}\frac{1}{n}-log(N)\right)\approx .577 215\ldots
\end{equation}

We do not know, even today (winter 2015), whether $\gamma$ is rational or not.\\

Later, in 1668 Mercator ( = Kremer) proved ([1], p. 185) the convergence of the alternative series (of even/odd numbers) and summed it (i.e., taking the $N\longrightarrow\infty$ limit):

\begin{equation}
1-\frac{1}{2}+\frac{1}{3}-\frac{1}{4}+\ldots+\frac{1}{2N-1}-\frac{1}{2N}+\ldots\longrightarrow\int_0^1\frac{dx}{1+x}=log(2)
\end{equation}

	Apparently, some work form India preceded the (even later) so-called  Gregory-Leibniz formula ([1] again, p. 184) for another alternative series (of inverse of odd numbers)

\begin{equation}
1-\frac{1}{3}+\frac{1}{5}-\frac{1}{7}+\ldots+\frac{1}{4N-3}-\frac{1}{4N-1}+\ldots\longrightarrow\int_0^1\frac{dx}{1+x^2}=\frac{\pi}{4}
\end{equation}

In this paper we \emph{interpret} (1.4) and (1.5) as arising from the cyclotomic equation of roots of the unity (Gauss; see, e.g. [5]):

 \begin{equation}
 \begin{aligned}
 &x^m-1=0\quad m\in\mathbb{N};\quad m=2,3,4,\ldots;\\
 &x^m-1= (x-1)(1+x+x^2 +...+ x^{m-1})\\
 &x^{2m}-1= (x-1)\left(1+x^2+x^4 +...+ x^{2(m-1)}\right)
\end{aligned}
\end{equation}
 	 	
	In our cases, $m = 2, 4$ respectively for (1.4) and (1.5), and then perhaps we can write anew the results in terms of some natural arithmetic functions. In (1.6), the roots   $\neq \pm1$ are complex conjugate pairs, as (1.6) describes a real equation, so the roots have modulus 1, and all of them lie in the complex unit circle $\approx S^1$.\\



	The main purpose of this paper is to generalize these results for other (generic) natural numbers $m\in\mathbb{N}$. We shall try to relate the series to some definite integrals.\\

	We anticipate already here part of the workings for $m$ = 2 and 4.\\

	For $m=2$, $x^2-1 = (x-1)(x+1)$; we \emph{always} leave out the $x=1$ root. The inverse of $(x+1)$ enters in (1.4) and we proceed to the four operations

\begin{description}
  \item[a)]\emph{inversion}, $(x+1) = (1+x) \longrightarrow(1+x)^{-1}$.	
  \item[b)] \emph{expansion}, $(1+x)^{-1}=1-x+x^2-x^3+x^4 -\ldots$  ($|x| < 1$).
  \item[c)] term-by-term \emph{integration}, $\longrightarrow x -\frac{x^2}{2} +\frac{x^3}{3} -\frac{x^4}{4} + \ldots$
  \item[d)] taking the \emph{limits} $x=1$ minus $x = 0$, $\longrightarrow(1- \frac{1}{2}) + (\frac{1}{3}-\frac{1}{4}) +(\frac{1}{5} -\frac{1}{6})\ldots$
\end{description}

Then we obtain the above result for $m = 2$:

\begin{equation}
\begin{array}{lcll}
&(x+1)=(1+x)\Longrightarrow\frac{1}{1+x}\Longrightarrow(1-x+x^2-\ldots)\Longrightarrow(x-\frac{x^2}{2}) + (\frac{x^3}{3}-\frac{x^4}{4}) +\ldots\Longrightarrow\\
&(\textrm{$x=1$ minus $x=0$})=\sum_{n=1}^\infty\frac{(-1)^{n+1}}{n}=log(2)
\end{array}
\end{equation}

\noindent in terms of the arithmetic modulated ``sign'' function (=even/odd) $f(n): =\frac{(-1)^{n+1}}{n}$. Note that \emph{after} the integration, the limit $x=1$ can be taken. (1.7) involves the series

\begin{equation}
\sum_{n=1}^{\infty}\left(\frac{1}{2n-1}-\frac{1}{2n}\right)=\int_0^1 dx \frac{1}{1+x}=\log(2)
\end{equation}

\noindent expressed as a genuine definite integral.\\

Similarly, for   $m = 4$, it is  $x^4-1= (x^2-1)(x^2+1)$, separating the roots $\pm 1$ and $\pm i$; now $(1+x^2)^{-1}$ enters in (1.5), which becomes, again after inversion, expansion, integration and taking limits:

\begin{equation}
(1+x^2)^{-1}\Longrightarrow (1-\frac{1}{3}) + (\frac{1}{5}-\frac{1}{7}) + (\frac{1}{9}-\frac{1}{11}) + \ldots ( = \frac{\pi}{4})
\end{equation}

\noindent which can be written as the (alternative) difference between the two series, and of course also as one integral

\begin{equation}
\sum_{n=1}^\infty\left(\frac{1}{4n-3}-\frac{1}{4n-1}\right)=\int_0^1 dx\frac{1}{1+x^2}=\frac{\pi}{4}
\end{equation}

And can be given now in terms of the so-called Dirichlet characters (see below).\\

	In this paper, as said, we shall generalize the constructions above (for $m = 2, 4$) for any integer $m\in\mathbb{N}$ (in principle), and include (when possible) the appropriate Dirichlet character.\\

	But first we want to fix our notation. For each $m\in \mathbb{N}$, $m >1$, we shall consider first the \emph{finite} sums, up to $mN$ terms; so we define primarily

\begin{equation}
\Sigma_m(N)\equiv 1 + \frac{1}{2} +\ldots+\frac{1}{mN-(m-1)}+\frac{1}{mN-(m-2)}  + \ldots + \frac{1}{mN}
\end{equation}
	
For example, $\Sigma_3(N)\equiv 1 + \frac{1}{2}+\frac{1}{3}+\ldots+\frac{1}{N}+\ldots +\frac{1}{3N-2}+\frac{1}{3N-1}+\frac{1}{3N}$; so these sums \emph{diverge}, for $N$ arbitrarily large, and indeed as (1.2) indicates:

\begin{equation}
\Sigma_m(N)\longrightarrow log(mN)+\gamma+\mathcal{O}(\frac{1}{mN})=log(N)+\gamma+log(m)
\end{equation}

\noindent (from now on we shall keep \emph{only} the \emph{dominant} and the \emph{constant} terms in the \emph{divergent} expansion  (1.2)).\\

	Also, we define \emph{partial bounded sums} $\Sigma^i_m=\Sigma^i_m(N)$ for $i= 1, 2, \ldots, m$) as displaced sums:

\begin{equation}
\begin{aligned}
&\Sigma^1_m = \frac{1}{1} + \frac{1}{m+1}  + \frac{1}{2m+1}+ \ldots + \frac{1}{(mN-(m-1))}\\	
&\Sigma^2_m = \frac{1}{2} + \frac{1}{m+2}  + \ldots + \frac{1}{(mN-(m-2))}\\	
\end{aligned}
\end{equation}
	   and\quad\quad\quad\quad\quad\quad\quad\quad$\ldots$

\begin{equation}
\Sigma_m^{m-1}=\frac{1}{m-1}+\frac{1}{2m-1}+\ldots+\frac{1}{mN-1}
\end{equation}

	And finally, $\Sigma_m^{m}=\frac{1}{m}+\frac{1}{2m}+\ldots+\frac{1}{mN}= \frac{1}{m}(1 + \frac{1}{2} + \ldots+ \frac{1}{N})$  so              	
\begin{equation}
(\textrm{for}\quad N\gg1)\quad \Sigma_m^{m}(N)\longrightarrow\frac{log(N)}{m}+\frac{\gamma}{m}+\mathcal{O}\left(\frac{1}{N}\right)
\end{equation}

Therefore, $\Sigma_m(N)$ (in 1.11) is the sum of \emph{all} $\Sigma_m^{i}(N)$ with $i : \{1, 2, \ldots, m\}:$\\

\begin{equation}
    \Sigma_m(N)= Sum_{i=1}^{m} \Sigma_m^{i}(N)=\Sigma_m^{1} +\Sigma_m^{2} +\ldots+\Sigma_m^{m}		                     \end{equation}

Note these sums make sense, spite diverging, because all have only positive terms and they are finite sums\ldots\\

The convergence or divergence of the above series is usually self-explained: convergence occurs \emph{always} for the decreasing alternative series (Leibniz), but as the convergence is not absolute, but conditional, the ordering in the series should be maintained. The higher divergence behaves, if at all, with a factor $\propto log(N)$: divergences are no higher than logarithmic. We shall try to be careful in subtracting two divergent series... We use a philosophy close to physics (in particular, to quantum electrodynamics or q. e. d.): we first truncate the series, taking a \emph{fixed} upper bound $N\gg 1$ (called the \emph{cut-off}; in physics this process is called \emph{regularization}).  Then we subtract other series, also divergent but with a similar type of divergence (this is called \emph{renormalization}): the result should be convergent (\emph{radiative} corrections); see e.g. the book by J. Schwinger [6].\\

	For a modern treatment of the Dirichlet series consult [7].\\

\section{The case for $m=2$: Four methods}.

 Here we repeat the $log(2)$ result for the case $m = 2$ (1.7) by \emph{four} different methods, because eventually the four might be useful.\\

	1) \emph{Redundance (direct) proof}: $\sum_2(N)$ and $\sum_2^2(N)$, defined above, diverge in a known way, for fixed (and large) $N$ and, as said (1.11, 1.14)

\begin{equation}
\begin{aligned}
&\Sigma_2^1=\Sigma_2-\Sigma_2^2\longrightarrow log(N) + log(2) + \gamma-\frac{1}{2}log(N)-\frac{\gamma}{2}=\\
&=\frac{1}{2}log(N)+\frac{\gamma}{2}+log(2)
\end{aligned}
\end{equation}

So
\begin{equation}
\Sigma_2^1-\Sigma_2^2\longrightarrow (N\longrightarrow\infty)=log(2)
\end{equation}

directly, without any integration or  (infinite) series summation.\\

	2) \emph{Integration proof}: Trivial here, as it is already (1.4):

\begin{equation}
\int_0^1\frac{dx}{1+x}=log(2)
\end{equation}

	3) \emph{Summation proof}: we have the \emph{convergent} series (so now $N =\infty$)

\begin{equation}
\begin{aligned}
&\Sigma_2^1-\Sigma_2^2=\left(1-\frac{1}{2}\right)+\left(\frac{1}{3}-\frac{1}{4}\right)+\ldots=\sum_{n=1}^{\infty}\left(\frac{1}{2n-1}-\frac{1}{2n}\right)=\\
&=\frac{1}{4}\sum_{n=1}^{\infty}\left(\frac{1}{n(n-\frac{1}{2})}\right)
\end{aligned}
\end{equation}

\noindent which can be explicitly done using the \emph{digamma function} $\Psi(z)$  [see e.g. [8], pag. 258]:

\begin{equation}
\Psi(z):=\frac{d(log\Gamma(z))}{dz}=\frac{\Gamma'(z)}{\Gamma(z)}\quad (z\notin\mathbb{Z}^-:=\mathbb{Z}\backslash\mathbb{N})
\end{equation}

\noindent (for $\Gamma(z)$ is the ordinary, Euler's gamma function), with help of the expression [8, p. 259, formula 6.3.16]

 \begin{equation*}
 \sum_{n=1}^{\infty}\frac{z}{n(n+z)}=\Psi(z)+\gamma+\frac{1}{z}
\end{equation*}

\noindent for $z = -\frac{1}{2}$ . In total we get, as expected

\begin{equation}
 \sum_{n=1}^{\infty} \frac{1}{n(n-\frac{1}{2})} = 4\log(2)
\end{equation}

4) Use of \emph{Hansen formula}: This formula, again depending on the Digamma function $\Psi(z)$, reads (Cfr. [9]; notice the sum starts in zero; $x$, $y$, $z$ are some parameters)

 \begin{equation}
 \sum_{n=0}^{\infty}\frac{1}{(nx+y)^2-z^2} =\frac{1}{2xz}\left\{\Psi\left(\frac{y+z}{x}\right)-\Psi\left(\frac{y-z}{x}\right)\right\}
 \end{equation}

We put our $n'$ above in $\sum_1\frac{1}{n'(n'-\frac{1}{2})}$ as ($n'= n+1$), because $n$ runs now from zero; hence this is equivalent to (2.7) with the parameter values $x=1$, $y=\frac{3}{4}$ and $z =\frac{1}{4}$. So our sum is (with $\Psi(1)=-\gamma$ and $\Psi(\frac{1}{2})=-\gamma-2\log(2))$

 \begin{equation}
 \Sigma= 2\left[\Psi(1)-\Psi\left(\frac{1}{2}\right)\right]=4\log(2)		
\end{equation}

\noindent as it should be.\\

Of the four methods, the most common and ``easy'' is the second (integration); it can, in principle (that is, if the integrand is known and the integration feasible), always be used. The \emph{Mathematica} Computer Programs give many integrals and double sums also directly.\\

	Please note in this (initial) $m = 2$ case that what appears as result is the \emph{sign function}$f(n):=\sum_{n=1}^{\infty}\frac{(-1)^{n+1}}{n}$; later, for some $m >2$, (including $m=4$) this will become a Dirichlet character: Here we have, as r\'{e}sum\'{e}  (to repeat)

 \begin{equation}
 \Sigma_2^1 -\Sigma_2^2 ( = log(2)) = \sum_{n=1}^{\infty}\frac{(-1)^{n+1}}{n}\equiv\Sigma^{\textrm{odd}}-\Sigma^{\textrm{even}}
\end{equation}

\noindent so NO Dirichlet character this time.\\

 \section{The case for  m = 3, 4, 6}.

	We discuss $m=3$ first. Now $(x^3-1) =(x-1)(1+x+x^2)$ as the roots  are +1, $\omega= \exp{\frac{2\pi i}{3}}$ and conjugate $\bar{\omega}$; here $\omega$ is equivalent to a plane rotation by $120^{\circ}$. By direct integration, we obtain at once

\begin{equation}
\int_0^1\frac{dx}{1+x+x^2}=\frac{\pi}{3\sqrt{3}}
\end{equation}

\noindent which is a particular case of

\begin{equation*}
\int \frac{1}{A}dx\quad\textrm{for}\quad A=(x-\omega)(x-\bar{\omega})=x^2-(Tr\omega)x+1,\quad Tr\omega=\omega+\bar{\omega}
\end{equation*}

For a general \emph{positive} integer $m$ we obtain, after an easy calculation ($\omega$ is root of $x^m-1=0$)

\begin{equation}
\int_0^1\frac{dx}{1-Tr\omega x+x^2}=\frac{\pi}{2 sin(\frac{2\pi}{m})}\cdot\frac{m-2}{m}
\end{equation}

This formula (3.2) can be used at once for $m =3,4$, and 6, yielding (3.1) and	

\begin{equation}
\int_0^1\frac{dx}{1-Tr\omega\cdot x+x^2}=\frac{\pi}{4}(m=4)\quad\textrm{and =}\frac{2\pi}{3\sqrt{3}}\quad (m=6).
\end{equation}

Later, we shall need also $I=\int_0^1 dx\frac{1}{1-Tr\omega x+ x^2}$ for $\omega =\exp{\frac{2\pi i r}{m}}$, for $r$ natural between 1 and $m-1$. The integral is now

\begin{equation}
I=\frac{\pi}{(2 sin(\frac{2\pi r}{m}))}\cdot\frac{(m-2r)}{m}\quad( 1\leq r \leq (m-1))	
\end{equation}

Coming back to $m = 3$, the specific series is  $1 -S + S^2 - S^3 + S^4 -\ldots$, with $S\equiv x+x^2$; developing, we get (after term-by-term integration, and taking limits 1, 0)

\begin{equation}
\begin{aligned}
&(1-\frac{1}{2}) + (\frac{1}{4}-\frac{1}{5}) + \ldots + \frac{1}{3N-2}-\frac{1}{3N-1} + \ldots =\\
&\sum_{n=1}^{\infty} \left(\frac{1}{3n-2}-\frac{1}{3n-1}\right)\equiv\Sigma_3^1-\Sigma_3^2=\frac{\pi}{3\sqrt{3}}
\end{aligned}
\end{equation}

Hence $\sum_3^1-\sum_3^2$ finishes the calculation for   $m = 3$,  (as the third sum, $\sum_3^3$ trivially computable).\\

In terms of the [10] Dirichlet character mod 3, namely $\chi^{(3)}_2$, we have

\begin{equation*}
\Sigma_3^1 - \Sigma_3^2 = \frac{\pi}{3\sqrt{3}} = \sum_{n=1}^{\infty}\frac{\chi^{(3)}_2 (n)}{n}
\end{equation*}

\noindent where, of course $\chi^{(3)}_2(1,2,3,4,5,6)=(1,-1,0,1,-1,0)$ (so 3-periodic). We can even use the Dirichlet's L-functions, which will include the $\frac{1}{n}$ factors, but we refrain of doing that, as it does not illuminate the matter any further.\\

As \underline{recapitulation}, the series are obtained from the polynomial of roots $\neq 1$

\begin{equation}
(m=3)\quad\quad P_2(x)=x^2+x+1
\end{equation}

\noindent after INVersion, EXPansion, INTegration and TAKing $x=1$. In principle, the result of the series summation can be also obtained from the Hansen formula (2.7).\\

Finally for this $m = 3$ case, we quote the four \emph{diverging} series for completeness:

\begin{equation}
		\Sigma_3^1-\Sigma_3^2 =\sum_{n=1}^{\infty} (\frac{1}{3n-2}-\frac{1}{3n-1})= \frac{\pi}{3\sqrt{3}} = \sum_{n=1}^{\infty}\frac{\chi^{(3)}_2 (n)}{n}\\
\end{equation}
as the \emph{complete} solution for the $m=3$ case. We quote the following four \emph{divergent} series for later use, still in this $m=3$ case (the limit $\longrightarrow$ meaning just  $N >>1$):

\begin{equation}
\Sigma_3 = 1 + \frac{1}{2} + \ldots + \frac{1}{3N} \longrightarrow log(N) + \gamma + log(3)
\end{equation}

\begin{equation}
\Sigma_3^3 =\frac{1}{3}+\frac{1}{6} + \ldots + \frac{1}{3N} \longrightarrow \frac{1}{3}log(N) + \frac{\gamma}{3}
\end{equation}

\begin{equation}
\Sigma_3^1 =1+\frac{1}{4}+\frac{1}{7} + \ldots + \frac{1}{3N-2} \longrightarrow \frac{1}{3}log(N) + \frac{\gamma}{3}+\frac{1}{2}log(3) + \frac{\pi}{6\sqrt{3}}
\end{equation}

\begin{equation}
\Sigma_3^2 =\frac{1}{2}+\frac{1}{5} + \ldots + \frac{1}{3N-1} \longrightarrow \frac{1}{3}log(N) + \frac{\gamma}{3}+\frac{1}{2}log(3) -\frac{\pi}{6\sqrt{3}}
\end{equation}\\

Now, for $m=4$, we have first the natural cyclotomic expression

\begin{equation}
			x^4-1 = (x^2-1)(x^2+1)\quad\quad\quad                I		              	
\end{equation}
(we write $I$ because this is \emph{not} the only possible factorization to be used). Repeating the steps as before for $m=3$, our first final result is here

\begin{equation}
\begin{aligned}
&\Sigma_4^1-\Sigma_4^3 =\sum_{n=1}^{\infty} \left(\frac{1}{4n-3}-\frac{1}{4n-1}\right)=\int_0^1\frac{dx}{x^2+1}=\\
&=arctan(1)=\frac{\pi}{4}
\end{aligned}
\end{equation}

with $\chi_2^{(m=4)}(1, 2, 3, 4, 5, 6, 7, 8) = (1, 0, -1, 0, 1, 0, -1, 0,)$ is periodic mod 4 (and restricted multiplicative). Note also why do we get $i = 1$ and 3 in $\Sigma_4^i$  (not 2!) mod 4: the expansion is for $\frac{1}{1+x^2}\approx 1-x^2+x^4-x^6$ etc., so it is with \emph{even} powers only, so with only \emph{odd} powers after integration!\\

We are done with this, as $\Sigma_4$, $\Sigma_4^4$ and $\Sigma_4^2$ are automatically obtainable.\\

The \underline{second} factorization of $x^4-1$ is obtained by separating only the $x=1$ root:

\begin{equation}
x^4-1= (x-1)(1+x+x^2+x^3)\quad\quad\quad		II
\end{equation}

As $B:= (1+x+x^2+x^3)$ contains the $x= -1$ root, one writes $B = (x+1)(1+x^2)$, where the integral can be computed at once (indeed, it is indicated above). The final result for factorization $II$ is

\begin{equation}
		\int_0^1 \frac{dx}{1+x+x^2+x^3} = \frac{1}{4}log(2) + \frac{\pi}{8}
\end{equation}

\noindent which is a \emph{redundant} result, as defines $\Sigma_4^1-\Sigma_4^2$ computable from the above calculation I. In fact, we end up this $m=4$ case by writing the analogous to (3.8) to (3.11):

\begin{equation}
\Sigma_4^4 =\frac{1}{4}+\frac{1}{8} + \ldots + \frac{1}{4N} \longrightarrow \frac{1}{4}log(N) + \frac{\gamma}{4}
\end{equation}

\begin{equation}
\Sigma_4^2 =\frac{1}{2}+\frac{1}{6} + \ldots + \frac{1}{4N-2} \longrightarrow\frac{1}{2}\Sigma_2^1= \frac{1}{4}log(N) + \frac{\gamma}{4}+\frac{1}{2}log(2)
\end{equation}

\begin{equation}
\begin{aligned}
&\Sigma_4^1 =\frac{1}{2}+\frac{1}{5} + \ldots + \frac{1}{3N-1} \longrightarrow\\
&\longrightarrow \frac{1}{4}log(N) + \frac{\gamma}{4}+\frac{3}{4}log(2) +\frac{\pi}{8}
\end{aligned}
\end{equation}

\begin{equation}
\Sigma_4^3 \longrightarrow \frac{1}{4}log(N) + \frac{\gamma}{4}+\frac{3}{4}log(2) - \frac{\pi}{8}
\end{equation}

The non-prime structure of $m$ \emph{implies} only a difference calculation, as it was also the case for $m=3$, prime. We have $\phi(3)=2$, $\phi(4)=2$.\\

So the full solution for the $m=4$ case has one redundancy (two factorizations), and it is done (solved) once a single calculation is made, e.g. $\int_0^1 dx \frac{1}{1+x^2} = \pi/4$.\\

Now it is the turn of \underline{$m=6$}. But here there are also several factorizations, as for the $m=4$ case above. The simplest is (perhaps)

\begin{equation}
\textrm{Firstly is}\quad x^6-1= (x^3-1)(x^3+1)\quad\quad       			 I
\end{equation}

As $(1+x^3) = (1+x)(1-x+x^2)$, the integral is easy, with this factoring:

\begin{equation}
\int_0^1\frac{dx}{1+x^3}=\frac{1}{3}log(2) + \frac{\pi}{3\sqrt{3}}
\end{equation}

The operations INV, EXP, INT and TAK $x=1$ applied to the polynomial $P_3(x) = 1+x^3$ yield the infinite but convergent sum

\begin{equation}
		\Sigma_6^1-\Sigma_6^4 =\sum_{n=1}^{\infty} \left(\frac{1}{6n-5}-\frac{1}{6n-2}\right)=\frac{log(2)}{3}+\frac{\pi}{3\sqrt{3}}
\end{equation}

Notice again the \emph{jump}, now by three: is due to the cubic $x^{3n}$ terms in $\frac{1}{1+x^3}$.\\

All other factorizations are therefore \emph{redundant}; we just write them:

\begin{equation}
(x^6-1) = (x-1)(1+ x+x^2+x^3+x^4+x^5) \Longrightarrow \Sigma_6^1-\Sigma_6^2\quad\quad 		I I
\end{equation}

\begin{equation}
(x^6-1) = (x^2-1)(1+x^2 +x^4)\Longrightarrow \Sigma_6^1-\Sigma_6^3\quad\quad   		        I I I
\end{equation}

Redundancy arises because from (3.22) one obtains $\Sigma_6^1$, as $\Sigma_6^4$ is directly computable (our Sect. 2); so, in $I I$, we know already the result, as $\Sigma_6^2$ is again trivial; same in $I I I$, as $\Sigma_6^3$ is again directly computable.\\

It is remarkable here that in the original expansion $(x^6-1) =\newline(x^3-1)(x^3+1)$ the simplest factorization $(x^3+1) =(x+1)(1-x+x^2)$, yields for $P_2(x) = 1-x+x^2$, the double series expansion

\begin{equation}
\begin{aligned}
&(1+\frac{1}{2})-(\frac{1}{4}+\frac{1}{5})+(\frac{1}{7}+\frac{1}{8})-(\frac{1}{10}+\frac{1}{11})\pm\ldots=\Sigma_6^1+\Sigma_6^2-\Sigma_6^4-\Sigma_6^5=\\
&=\int_0^1\frac{dx}{1-x+x^2}=\frac{2\pi}{3\sqrt{3}}
\end{aligned}
\end{equation}

\noindent which again it is a \emph{redundant} calculation, as included in the (summed)  $\int_0^1 dx \frac{1}{1+x^3}$.\\

To finish these simple cases, we just write down the nontrivial summations:

\begin{equation}
\begin{aligned}
&\Sigma_6^1(N)  = \frac{1}{6}log(N) + \gamma/6+\frac{log(2)}{3} + \frac{log(3)}{4}+\frac{log(2)}{3} + \frac{\pi}{4\sqrt{3}}\\
&\Sigma_6^4 = \frac{log(N)}{6} +\frac{\gamma}{6} +\frac{log(3)}{4}-\frac{\pi}{12\sqrt{3}}\\
&\Sigma_6^1-\Sigma_6^5=\frac{log(3)}{4}
\end{aligned}
\end{equation}

And, in this  $m =6$ case, we get also\\

$\Sigma_6^1-\Sigma_6^4=\frac{1}{3}\log(3)+\frac{\pi }{3\sqrt{3}}$.\\

\noindent computing $\Sigma_6^2$ and $\Sigma_6^4$ as in  Sect. 1.\\

	In terms of some Dirichlet functions we get\\
\begin{equation}
\Sigma_6^1-\Sigma_6^5=\Sigma\frac{\chi^{(6)}_2}{n},\quad \textrm{as}\quad \chi^{(6)}_2(1,2,3,4,5,6)=(1,0,0,0,-1,0)=\frac{log(3)}{4}
\end{equation}

For a modern treatment of sums involving harmonic numbers, see [4].\\

\section{General numbers}.

Here is the cyclotomic equation is, for $m=5$,

\begin{equation}
\begin{aligned}
&x^5-1=(x-1)(1+x+x^2+x^3+x^4)=\\
&=(x-1)(x-\omega)(x-\bar{\omega})(x-\omega^2)(x-\bar{\omega^2})\\
\end{aligned}
\end{equation}

\noindent where $\omega= \exp{\frac{2\pi i}{5}}$ = Rotation by $72^{\circ}$ and $\omega^2 =\exp{\frac{4\pi i}{5}}$ = Rotation by $144^{\circ}$.\\

We obtain easily
			
\begin{equation}
cos(72^{\circ}) = \frac{\sqrt{5}-1}{4}\quad  cos(144^{\circ}) =\frac{-1-\sqrt{5}}{4}
\end{equation}

The integral $T \equiv \int_0^1 \frac{dx}{1+x+x^2+x^3+x^4}$ can be done as the real denominator splits in two real ones and quadratic, but we abstain to write explicitly $(\omega=e^{\frac{2\pi i}{5}})$

\begin{equation}
(1+x^2+x^3+x^4)=(1-xTr\omega +x^2)(1-xTr\omega^2 +x^2)
\end{equation}

The expression for the integral is too long to write it. In summation terms, it is
$\frac{1}{1+x+x^2+x^3+x^4} = \ldots = (1-\frac{1}{2}) + (\frac{1}{6}-\frac{1}{7}) + \ldots =\sum_{n=1}^{\infty} \left(\frac{1}{5n-4}-\frac{1}{5n-3}\right)$ still computable (e.g. with Mathematica), but still too long. The final result will be

\begin{equation}
\int_0^1dx\frac{1}{1+x+x^2+x^3+x^4}=\Sigma_5^1-\Sigma_5^2
\end{equation}

This does NOT correspond to none mod 5 Dirichlet characters, but one can always put a (Dirichlet) function

\begin{equation*}
f^{(5)}_2(1,2,3,4,5,6,7\ldots)=(1,-1,0,0,0,1,-1)\quad \textrm{periodic mod 5}
\end{equation*}

The simplest summation (according to \emph{Mathematica}) is

\begin{equation}
\Sigma_5^1-\Sigma_5^4=\frac{(1+\sqrt{5})\pi}{5\sqrt{2(5-\sqrt{5})}}
\end{equation}

Notice the factorization in (4.3) does not allow us to write it as difference between two series, because the coefficients are not integer numbers.\\

For the next prime, namely $m =7$, we have three couples of complex roots, plus the $x=1$ value: the sextic integral has not been attempted, but the summation can be again being done; we refrain to elaborate.\\

This is the general trend for prime numbers $p$; there are $\frac{(p-1)}{2}$ pairs of complex conjugate roots; ``a priori'', the only integral vs. series is the simplest case, generalizing the above result

\begin{equation}
\int_0^1dx\frac{1}{1+x+x^2+\ldots+x^{p-1}}=\Sigma_p^1-\Sigma_p^2
\end{equation}

\section{General case}.

We deal now with some composite numbers.
Compound numbers are easier; we just add a calculation for $m = 8$, with $\phi(8) = 8-4 = 4$:\\

\begin{equation}
x^8-1= (x^4-1)(x^4+1),\quad  \textrm{and}\quad (x^4+1) = (x^2-\sqrt{2}x +1)(x^2 + \sqrt{2}x +1),\quad \textrm{and so}
\end{equation}

\begin{equation}
\begin{aligned}
\int_0^1\frac{dx}{x^4+1}&=\frac{\pi+\log{(3+2\sqrt{2})}}{4\sqrt{2}}=\ldots=\Sigma_8^1-\Sigma_8^5=\\
&=\sum_{n=1}^{\infty} \left(\frac{1}{8n-7}-\frac{1}{8n-3}\right)=\frac{\pi}{4\sqrt{2}}+\frac{\log(3+2\sqrt{2})}{4\sqrt{2}}
\end{aligned}
\end{equation}

Another factorization is $(x^8-1) = (x^2-1)(1+x^2+x^4+x^6)$. It is equivalent to\\

\begin{equation}\int_0^1 dx\frac{1}{(1+x^2+x^4 + x^6)} = \Sigma_8^1-\Sigma_8^3 = \sum_{n=1}^{\infty} \left(\frac{1}{8n-7}-\frac{1}{8n-5}\right)
\end{equation}

Still, a third factorization is $(x^8-1) =(x-1)(1+x+x^2+\ldots + x^7)$, which yields $\Sigma_8^1-\Sigma_8^2$, so the remaining $\Sigma_8^7$ is obtained by difference (with $\Sigma_8^{2,4,6,8}$ inmediate). The case $m=8$ is potentially resolved. Again, we refrain to elaborate.\\

For $m = 9$, $\Sigma_9^3$ and $\Sigma_9^6$  can be computed directly, whereas the odd cases $\Sigma_9^{1,3,5,7}$ require further work, but it is again feasible. Similarly for $m =10.$\\

As \emph{general conclusion}, we have showed a remarkable relation between the cyclotomic equation $x^m-1=0$ and some series and definite integrals; they go from the simplest integrals (and series) in the literature (like $\int dx \frac{1}{1+x}=\log(2)$) to very complicated cases, still feasible: the integrals have denominators factoring in quadratic ones, and the series are of the type $\sum \frac{1}{an^2 + bn + c}$, computable, in principle, by means of the Hansen's formula.\\

	There are, however, some questions left over in our work: for example, the series for   $\frac{1}{1+x+x^2+x^3+\ldots+x^q}$ we identify with the series $\sum_{q+1}^1$; this is correct, but we have checked it ``case by case'', offering no general proof; etc. Also we feel some new series might perhaps appear, whenever the quadratic components offer an integer expansion: those are two questions for the future...

\section{ACKNOWLEDGEMENTS.}

. We are very thankful to Prof. Jes\'us Guillera, from the Mathematics Department for his interest and for showing the authors the relations with the Dirichlet characters and also the use of the Hansen formula. This work has been supported by the Spanish CICYT (grant FPA 2006-02315) and the aragonese DGIID-DGA (grant 2007-E242).

\section{REFERENCES}

[1].- J. Stillwell, \emph{Mathematics and Its History}. Springer 2010.\\

[2].- L. Jolley, \emph{Summation of series}. Dover 1961.\\

[3].- T. J. Bromwich, \emph{Introduction to The theory of Infinite Series}.\\
\quad\quad\quad        McMillan 1926 				       (courtesy of University of Toronto Library).\\

[4].- W. Chu, \emph{Summation formulae involving harmonic numbers}. Filomat 26:1 (2012), 143-152.\\

[5].-H. Kragh S{\o}rensen \emph{Gauss and the cyclotomic equation}. Preprint (P.h. D. Thesis), Aarhus (Denmark) University, 3-III-2001.\\

[6].- J. Schwinger (ed.) \emph{Selected papers in Quantum Electrodynamics}. Dover 1956.\\

[7].- R. J. Mathar, \emph{Table of Dirichlet L-Series and Prime Zeta modulo \quad\quad\quad\quad\quad\quad\quad\quad\quad\quad \quad\quad\quad\quad\quad  Functions for Small Moduli}. arXiv: 1008.2547 [math. NT], 15-VIII-2010.\\

[8].- M. Abramowitz and I. Stegun, \emph{Handbook of Mathematical Functions}. Dover 1965.\\

[9].- E. R. Hansen, \emph{A table of series and products}, p. 105. Prentice-Hall, 1975.\\

[10].- T. Apostol, \emph{Introduction to Analytic Number Theory}.\\
\quad\quad Springer 1976.
\end{document}